\documentclass[12pt]{article}
\usepackage{amssymb,amsmath,amsthm}

\theoremstyle{plain}
\newtheorem{theorem}{Theorem}[section]

\newtheorem{lemma}[theorem]{Lemma}

\theoremstyle{definition}
\newtheorem{defn}[theorem]{Definition}
\newtheorem{remark}[theorem]{Remark}

\newcommand{\Hn}{\mathbb{H}^n} 
\newcommand{\Hk}{\mathbb{H}^k}
\newcommand{\vol}{\mathrm{vol}} 
\newcommand{\bary}{\mathrm{bar}}
\newcommand{\isom}{\mathrm{Isom}(\Hn)}
\newcommand{\e}{\ensuremath{\varepsilon}}
\newcommand{\B}{\ensuremath{\mathcal B}}

\title{Maximal volume representations are fuchsian}

\author{Stefano Francaviglia\footnote{The first author was supported
    by an INdAM fellowship and the Department of Applied Mathematics
    of the University of Pisa.}\\{\footnotesize University of Pisa, Italy}\\
    {\footnotesize e-mail: francavi@sns.it} \and Ben Klaff\footnote{The second
    author was also supported by a CIRGET fellowship and by the Chaire de Recherche du Canada en alg\`ebre, combinatoire et informatique math\'ematique de l'UQAM.}\\
    {\footnotesize University of Texas, Austin}\\{\footnotesize e-mail: klaff@math.utexas.edu }\\
}

\date{November 2, 2004}

\begin{document}

\maketitle

\begin{abstract}
We prove a volume-rigidity theorem for fuchsian representations
of fundamental groups of hyperbolic $k$-manifolds into $\isom$.  Namely,
we show that if $M$ is a complete hyperbolic $k$-manifold with finite volume,
then the volume of any representation of $\pi_1(M)$ into $\isom$, $3 \leq k \leq n$, is less than the volume of $M$, and the volume is maximal if and only if the representation is discrete, faithful and ``k-fuchsian''.
\end{abstract}

\section{Introduction} \label{intro}

The main result of this paper is a generalization and streamlined
proof of a result which is often referred to as the ``representation
volume rigidity'' theorem:

\begin{theorem} \label{mainthm}
Let $M$ be an oriented, connected, complete, real hyperbolic $k$-manifold of finite volume, with $k\geq 3$.
Let $\rho: \pi_1(M) \rightarrow \isom$ be a representation of its
fundamental group into the group of isometries
of hyperbolic $n$-space.  Then the volume of
$\rho$ is less than or equal to the volume of $M$, and equality holds
if and only if $\rho$ is $k$-fuchsian, i.e., a discrete and faithful
representation into the group of isometries of a $k$-dimensional
subspace of $\Hn$.
\end{theorem}

In the case $k=n=3$ and $M$ is closed, this result was proved, following original ideas of W. Goldman, M. Gromov, and W. Thurston~\cite{Gol82,Thu:note}, by N. Dunfield~\cite{Dun99}.  The result was extended to the case when $M$
has finite volume by the authors~\cite{Fra04, Kla:tesi}.

The new ingredient in the present proof is the use of natural maps, or the
barycenter method, a technique introduced and developed by G. Besson, G. Courtois, and S. Gallot.  (See, for example,~\cite{BCG95,BCG96,BCG99}.)  A key step in the proof is the first author's generalization of the B-C-G method
so as to be able to construct natural maps for representations~\cite{F5}.

We remark that in addition to the new proof offered here, the
original ``volume-rigidity of representations'' result itself has been generalized, in that it deals with the case when the target dimension is greater than that of the domain.  We also note that Besson, Courtois and Gallot have recently
obtained a similar result~\cite{BCG04}.

The paper is organized as follows.  In Section~\ref{defs}, we give the
necessary background and definitions, including the definition of
the volume of a representation.  In Section~\ref{proofmainthm}, we give
the proof of Theorem~\ref{mainthm} in the case $M$ is closed.  Finally, in Section~\ref{finitevolume}, we complete the proof of Theorem~\ref{mainthm},
giving the proof for non-compact, complete, finite-volume manifolds.

\vskip\baselineskip\noindent\textsc{Acknowledgements.}
The authors warmly thank Juan Souto for the stimulating conversations about
the B-C-G techniques while in Cambridge and in Pisa.  The second author would also like to thank Pete Storm, Benson Farb, Steve Boyer, and Stephan Tillmann, as well as the Universita di Pisa for their hospitality.

\section{Definitions and notation} \label{defs}

Throughout Sections~\ref{defs} and~\ref{proofmainthm}, $M$ will denote
a closed, oriented hyperbolic $k$-manifold, $k \geq 3$.  We
suppress a choice of basepoint in $M$ and let $\Gamma = \pi_1(M)$
denote the fundamental group of $M$.  We let $\rho: \Gamma \rightarrow
\isom$ denote a representation of $\Gamma$ into the group of
isometries of $\Hn$.

By a {\it pseudo-developing map} for $\rho$, we mean a piecewise
smooth map $D=D_{\rho}: \widetilde{M} \rightarrow \Hn$ which is
$\rho$-equivariant, i.e., such that
\[ D (\gamma \cdot x) = \rho(\gamma) \cdot D(x),\]
for every $x \in \widetilde{M} = \Hk$, the universal cover of $M$, and
for every $\gamma \in \Gamma$.

Given any representation $\rho$, one can construct a
pseudo-developing map for $\rho$ as follows: lift a smooth
triangulation for $M$ to $\widetilde{M}$, and then recursively define
the map $D$ on the $i$-skeleta, $0 \leq i \leq k$, by choosing images
for a complete system of orbit representatives for the $i$th skeleta,
and then extending the map equivariantly.

We now introduce the notion of the volume of a representation.  Let $h$ denote the hyperbolic metric on the target $\Hn$ and let $D$
be a pseudo-developing map for $\rho$.  The pullback $D^{*}h$ of $h$
along $D$ is a (possibly degenerate) pseudo-metric on $\widetilde{M}$ and
hence induces a $k$-form $\tilde{\omega}_{D} = |
\mathrm{det} \, D^{*}h|$ on $\widetilde{M}$.  Since the map $D$ is
$\rho$-equivariant, $D^{*}h$ and hence $\tilde{\omega}_D$ are
$\Gamma$-invariant.  Hence the $k$-form $\tilde{\omega}_D$ descends to
a $k$-form $\omega_D$ on $M$.

\begin{defn}[Volume of a pseudo-developing map]\label{pseudodevvoldefn}
The {\it volume} $\vol(D)$ {\it of a pseudo-developing map} $D$ for a
representation $\rho$ is defined by
\[ \vol(D) = \int_{M} \omega_D.\]
\end{defn}

We can now make the following

\begin{defn}[Volume of a representation]\label{repvoldefn}
The {\it volume $\vol(\rho)$ of a representation} $\rho$ is defined by
\[ \vol(\rho) = \mathrm{inf}_{D} \{ \vol(D) \},\]
where the infimum is taken over the set of all pseudo-developing maps
$D$ for $\rho$.
\end{defn}

Note that $\vol(D)$ and hence $\vol(\rho)$ are non-negative real
 numbers.  Also, note that $\vol(D)$ is not invariant under
 $\rho$-equivariant homotopy.  Hence the volumes of two
 pseudo-developing maps for a given representation can be different.
 We use the above definition of representation-volume
 in order to deal with the case $n \neq k$.  (Compare the
 definition of representation-volume and the consequent property of
 invariance under homotopy in~\cite{Dun99,Fra04}.)  Finally, we point out that in the non-compact case, the definition of volume of a representation involves another condition.  (See Section~\ref{finitevolume}.)

\section{The compact case}\label{proofmainthm}

When $M$ is compact, the proof of Theorem~\ref{mainthm} 
goes as follows.  First, we invoke an existence result due to the first
author~\cite{F5}, which says that there is a pseudo-developing map $F$ for
$\rho$ such that $\vol(F) \leq \vol(M)$.  The inequality then follows by
the definition of $\vol(\rho)$.

 Next, we use the hypothesis that $\vol(\rho) = \vol(M)$, some elementary
 Riemannian geometry, and the properties of the pseudo-developing map
 $F$ to conclude that $F$ is a Riemannian isometry from $\Hk$ to a
 $k$-dimensional hyperbolic subspace of $\Hn$.  (This then
 reduces the remainder of the proof to the case $k=n$.)  It is then easy to conclude that $F$ is a covering map onto its image, and it follows that $\rho$ is discrete and faithful. 

 Finally, we show the (easier) converse, namely that if $\rho$ is a discrete,
 faithful representation into the group of isometries of a
 $k$-dimensional hyperbolic subspace of $\Hn$, then $\vol(\rho) =
 \vol(M).$

The following result is proved in~\cite{F5}.

\begin{lemma}\label{mainlemma}
Let $\rho: \Gamma \rightarrow \mathrm{Isom}(\Hn)$ be a representation
whose image is a non-elementary group.  Then, for $k \geq 3$, there
exists a smooth pseudo-developing map $F : \Hk \rightarrow \Hn$ such
that for all $x \in \Hk$,
\begin{equation}\label{ineqone} |\mathrm{Jac} \,F(x) | \leq 1;\end{equation}
moreover, equality holds at $x$ if and only if $dF_{x}: T_{x}\Hk
\rightarrow T_{F(x)}\Hn$ is an isometry.
\end{lemma}  

Assuming the image of $\rho$ is non-elementary, Lemma~\ref{mainlemma}
now implies the inequality of Theorem~\ref{mainthm}.  Indeed, by the definition of volume of a representation and the
inequality in the lemma, it follows immediately that
\begin{equation}\label{ineqtwo} 
\vol(\rho) \leq \vol(F) \leq \vol(M).
\end{equation} 
 If the image of $\rho$ is elementary, then it is easy to check that $\vol(\rho) = 0$.  Thus, in either case, inequality (\ref{ineqone}) holds.

\vskip\baselineskip
 
We now suppose that $\vol(\rho) = \vol(M)$ and proceed to show
that the image of the pseudo-developing map $F$ is
contained in a $k$-dimensional hyperbolic subspace of $\Hn$.

Since $\vol(\rho) = \vol(M)$, each of the inequalities of
$(\ref{ineqtwo})$ is an equality.  Hence, for each $x$ in $\Hk$, the
inequality in $(\ref{ineqone})$ is equality.  Thus the map $F$ is a
Riemannian isometry.

We now recall some ideas and facts from Riemannian geometry, referring
the reader to~\cite{GHL:libro} for notation and details.  We note that, in what follows, $C^2$-regularity of the pseudo-developing map is enough.

Let $X$ denote a Riemannian manifold.  A submanifold $N$ of $X$ is called
{\it minimal} if it is a critical point of the volume function.  A
submanifold is {\it locally minimal} if, for each point $x$ of $N$,
there exists a neighborhood $A$ of $x$ such that all perturbations of
$N$ with support in $A$ do not decrease the volume of $N$.  A
submanifold $N$ of $X$ is {\it totally geodesic} if for any two points
$x$ and $y$ in $N$, the geodesic joining $x$ and $y$ in $X$ is
contained in $N$.  We denote by $R^{X}$ and $\nabla^{X}$ (resp., $R^N$
and $\nabla^N$) the curvature tensor and the connection of $X$ (resp.,
$N$).

For any two vector fields $U$ and $V$ in $N$, we denote by $\Pi(U,V)$
the second fundamental form of the submanifold $N$.  Equivalently, if
$\{ \nu_1\ \ldots, \nu_r\}$ denotes an orthonormal frame of the
orthogonal complement of $TN$ in $TX$, and if $l_{i}(U,V)$ denotes the
real-valued fundamental form corresponding to $\nu_{i}$, then
\[ \nabla_{U}^{X}V - \nabla_{U}^{N}V = - \sum_{i=1}^{r}l_{i}(U,V)\nu_{i} =
\Pi(U,V).\]

The strategy is to prove that the image of the map $F$ is a minimal
submanifold of $\Hn$, and from this conclude that the image of $F$ is
contained in a $k$-dimensional subspace of $\Hn$.  To do this, we need
the following standard results (\cite[Chapter V]{GHL:libro}).

\begin{lemma} \label{diffgeolemma1}
Let $N$ be a submanifold of a Riemannian manifold $X$.  Then
\begin{enumerate}
\item\label{minimalcrit1} 
$N$ is minimal if and only if the traces of all the real-valued
second fundamental forms vanish (see~\cite[p.228]{GHL:libro});
\item\label{minimalcrit2} $N$ is totally geodesic if and only if the second fundamental
form vanishes (see~\cite[p. 220]{GHL:libro}). 
\end{enumerate}
\end{lemma}

\begin{lemma} \label{diffgeolemma2}
The image $F(\Hk)$ of the map $F$ is contained in a locally minimal
submanifold of $\Hn$.
\end{lemma} 

\begin{proof}
Suppose not.  Then by a perturbation of $F$ in a small ball $B$ of
$\Hk$, we can decrease the volume of $F$.  Indeed, by
$\rho$-equivariantly perturbing $F$ in the $\Gamma$-orbit of $B$, we
can find a pseudo-developing map $F': \Hk \rightarrow \Hn$ with a
strictly smaller volume than that of $F$.  But then $\vol(M) =
\vol(\rho) \leq \vol(F') < \vol(F) = \vol(M)$, a contradiction.
\end{proof}

\begin{lemma} \label{diffgeolemma3}
Let $N$ be a locally minimal $k$-submanifold of a Riemannian
$(k+r)$-manifold $X$.  If, for all vector fields $U,V,W$, and $T$ we
have
\[ R^N(U,V,W,T) = R^{X}(U,V,W,T),\]
then $N$ is totally geodesic.
\end{lemma}

\begin{proof}
By $(\ref{minimalcrit2})$ of Lemma~\ref{diffgeolemma1}, it suffices to
 show that the second fundamental form of $N$ vanishes.  We again let
 $\{ \nu_1, \ldots, \nu_{r}\}$ be an orthonormal frame of the
 orthogonal complement $TN$ of $TX$, and for each index $i$, we let
 $l_{i}(\cdot,\cdot)$ denote the real-valued fundamental form
 corresponding to $\nu_{i}$.  By Gauss's theorem (see for
 example~\cite[Chapter V]{GHL:libro}), we
 conclude that for any point $p \in N$ and for any $u,v,w$, and $t$ in
 $T_{p}N$,
\[ R^{N}(u,v,w,t) = R^{X}(u,v,w,t) + \sum_{i=1}^{r} (l_{i}(u,w)l_{i}(v,t)
- l_{i}(u,t)l_{i}(v,w)).\] It then follows that for any $u,v,w$, and
$t$ in $T_{p}N$,
\[ \sum_{i=1}^{r}(l_{i}(u,w)l_{i}(v,t) - l_{i}(u,t)l_{i}(v,w))) = 0.\]
By hypothesis, $N$ is a locally minimal submanifold; therefore, by
$(\ref{minimalcrit1})$ of Lemma~\ref{diffgeolemma1}, we have that
$\textrm{tr}(l_{i}) = 0$ for each $i=1, \ldots, r$.

Now let $e_{1}, \ldots, e_{k}$ denote an orthonormal basis of
$T_{p}N$.  Setting $u=t=e_{j}$ in the above equality, we have
\[ \sum_{i=1}^{r}(l_{i}(e_{j},w)l_{i}(v,e_{j}) - l_{i}(e_{j},e_{j})l_{i}(v,w))) = 0.\]
Setting $w = v$ and summing over the index $j$, we get
\[ \sum_{j=1}^{k} \sum_{i=1}^{r} l_{i}^{2}(e_{j},w) - \sum_{i=1}^{r} \textrm{tr}(l_{i})l_{i}(w,w))
= 0.\] Whence, by the vanishing trace condition of
Lemma~\ref{diffgeolemma1}, we conclude that for any $p$ in $N$ and $w$
in $T_{p}N$,
\[ \sum_{i,j} l_{i}^{2}(e_{j}, w) = 0.\]
It now follows that $l_{i}(e_{j},w) = 0$ for any $i,j$, and $w$, and
hence that $l_i \equiv 0$ for $1 \leq i \leq r$.  This shows that the
second fundamental form vanishes at each point $p$ in $N$, which
completes the proof of the lemma.
\end{proof}

We now apply Lemma~\ref{diffgeolemma3} with $N = F(\Hk)$ and $X =
\Hn$.  Since $F$ is a Riemannian isometry, the hypothesis that $R^{N}
= R^{X}$ is satisfied.  By Lemma~\ref{diffgeolemma2}, $N$ is a locally
minimal submanifold of $X$.  Hence by Lemma~\ref{diffgeolemma3}, $N$
is totally geodesic.  Therefore the map $F$ is an isometry from $\Hk$
to a $k$-dimensional subspace $H$ of $\Hn$, and it follows that the
image of $\rho$ is contained in the group of isometries of $H$, as desired.

We claim now that $F: \Hk \rightarrow H$ is a covering map. Indeed, 
note that there exists an $r>0$ such that for any $x \in
\Hk$, the restriction map $F|_{B(x,r)}$ is an isometry onto its image.
This easily implies the claim.

Since $\Hk$ is simply connected, the covering $F:\Hk \rightarrow H$ is
a homeomorphism.  Thus $F$ is a $\rho$-equivariant global isometry of
$\Hk$.  It follows that the representation $\rho$ is discrete and
faithful.

Finally, we suppose that $\rho$ is a discrete and faithful
representation into the group of isometries of a $k$-dimensional
subspace $H$ of $\Hn$, and show that $\vol(\rho) = \vol(M)$.  First,
note that it is not restrictive to consider only those pseudo-developing maps for $\rho$ whose images are contained in $H$. 
Thus, after identifying $H$ with $\Hk$, we
may assume that $n=k$.

Let $N=\Hk / \rho(\Gamma)$.  By Mostow rigidity, the hyperbolic
$k$-manifolds $M$ and $N$ are isometric, and in particular, $\vol(N) =
\vol(M)$.  Now let $D$ be any pseudo-developing map for $\rho$.  Since
$D$ is $\rho$-equivariant, it induces a map $g: M \rightarrow N$, and
by definition, $\vol(D) = \int_{M} |g^{*}\omega|$, where $\omega$ is
the hyperbolic volume form of $N$.  Hence
\[ \vol(N) = \vol(M) = |\int_{M} g^{*}\omega| \leq \int_{M} |g^{*}\omega|
= \vol(D).\] It follows that
\[ \vol(M) \leq \vol(\rho) = \mathrm{inf}_{D} \{ \vol(D)\},\]
and we have already shown (see inequality $(\ref{ineqtwo})$ after Lemma
\ref{mainlemma}) that the reverse inequality also holds.  This
completes the proof of Theorem~\ref{mainthm} when $M$ is a compact manifold.

\section{The finite-volume case}
\label{finitevolume}
In this section we complete the proof of Theorem~\ref{mainthm},
proving the result in the finite-volume case. The main difference from the
previous case is that, as our manifolds are no longer
compact, we need to work with proper maps; since we work at the
level of universal coverings, we need an equivariant notion of
properness.  We keep here all the notation and definitions  of previous
sections, except that in the sequel $M$ will denote an oriented, complete,
non-compact, hyperbolic $k$-manifold of finite volume with $k\geq3$.  
We will also need to modify the definition of the volume of a
representation.

The manifold $M$ is diffeomorphic to the interior of a compact
manifold $\overline M$ whose boundary consists of Euclidean
$(k-1)$-manifolds.  (See, for example,~\cite{BePe:libro}).
In particular, for each boundary component $T\subset\partial \overline
M$ the group $\pi_1(T)<\pi_1(M)=\Gamma<\rm{Isom}(\Hk)$ is an abelian
parabolic group. The following lemma is easy to check.

\begin{lemma}
Let $G$ be an abelian group of isometries of a hyperbolic space
$\mathbb H^m$. Then the set $\rm{Fix}(G)\subset\overline{\mathbb H^m}$
of points  which are fixed by $G$ is non-empty.  
\end{lemma}

We note that $G$ may have no fixed point in $\partial \mathbb H^m$
(for example if $G<\rm{Isom}(\mathbb H^3)$ is the dihedral group
generated by two rotations of angle $\pi$ around 
orthogonal axes). 

Up to conjugacy, a peripheral subgroup of $\pi_1(M)$ has a unique fixed point, which lies in $\partial \Hk$.  Thus, for each $T\subset \partial\overline M$, each conjugate of $\pi_1(T)$ in $\pi_1(M)\subset\rm{Isom}(\Hk)$
corresponds to its fixed point in $\partial \Hk$.

We can now give the definition of a properly-ending map.
\begin{defn}[Properly ending maps]
  Let $\rho:\pi_1(M)\to\isom$ be a representation, and let $D:\Hk\to\Hn$
  be a $\rho$-equivariant map. We say that $D$ {\em properly ends} if
  for each $T\subset \partial \overline M$, if
  $\xi=\rm{Fix}(\pi_1(T))$ and $\alpha(t)$ is a geodesic ray ending at
  $\xi$, then all limit points of $D(\alpha(t))$ lie either in
  $\rm{Fix}(\rho(\pi_1(T)))\subset \overline\Hn$ or in a finite union
  of $\rho(\pi_1(T))$-invariant geodesics. 
\end{defn}

\begin{defn}[Volume of a representation]
The {\it volume $\vol(\rho)$ of a representation} $\rho$ is defined by
\[ \vol(\rho) = \mathrm{inf}_{D} \{ \vol(D) \},\]
where the infimum is taken over the set of all properly-ending
pseudo-developing maps 
$D$ for $\rho$.
\end{defn}

\begin{remark}
  It is easy to construct properly-ending pseudo-developing maps. (See~\cite{Dun99,Fra04}.)  We need to work with such maps because
  otherwise, one can construct (non-properly-ending) pseudo-developing maps with volume zero.  (For example, one can collapse $M$ to any of its spines.)  Also, we note that the above definition of volume
  ``extends'' the previous one given for compact manifolds.  Indeed, if
  $M$ is compact, then any pseudo-developing map properly ends.
\end{remark}

We now need to recall the definition and properties of the barycenter of
measures in $\overline\Hn$, referring to~\cite{BCG99,F5} for details.  (The reader who is familiar with such constructions may skip directly 
to Lemma~\ref{mainlemma2}.)  Let $\beta$ be a probability Borel measure
on $\partial \Hn$.  We define a function $\B_\beta:\Hn\to\mathbb R$ by
$$\B_\beta(y)=\int_{\partial \Hn}B(y,\theta)\, d\beta(\theta)\label{beta}$$
where $B(y,\theta)$ is the Busemann function of $\Hn$.
Then we have
\begin{enumerate}
\item If $\beta$ is not concentrated in two points, then $\B_\beta$ is
  strictly convex (because its Hessian is the $\beta$-average of the
  Hessians of the Busemann functions $B(y,\cdot)$) and goes to
  $\infty$ as $y$ goes to $\partial \Hn$.
\item If $\beta$ is not the sum of two Dirac delta measures with the same
  weight, then $\B_\beta$ has a unique minimum (possibly $-\infty$) in
  $\overline \Hn$.  Such a minimum is attained in $\partial \Hn$
  if and only if $\beta$ has an atom of weight greater that
  $\frac{1}{2}$.  The point $\bary(\beta)$ where $\B_\beta$ attains its
  minimum is called the {\em barycenter} of $\beta$.
\item\label{newpoint3} If $\beta$ is the sum
  $\frac{1}{2}(\delta_{\theta_1}+\delta_{\theta_2})$ of two Dirac
  delta measures concentrated in $\theta_1$ and $\theta_2$, then $\B_\beta$ is convex and constant on the geodesic joining $\theta_1$ and
  $\theta_2$, where it attains its minimum.
\item If $\beta$ is a probability measure on $\Hn$, its barycenter is
  defined by taking the convolution with the family of visual measures as
  follows.  Let $\nu_{O'}$ be the standard probability measure on
  $\partial \Hn\simeq \mathbb S^{n-1}$ in the disc model with center
  $O'$.  For every $y\in\Hn$, define $\nu_y=\psi_*\nu_{O'}$, where $\psi$
  is any isometry mapping $O'$ to $y$.  (Note that this is well-defined
  because $\nu_{O'}$ is Stab$(O')$-invariant.)  Now define $\bar\beta$,
  a probability measure on $\partial \Hn$, by 
$$\int_{\partial \Hn}\varphi(\theta)\,d\bar\beta(\theta)=
\int_{\Hn}\left(\int_{\partial \Hn}\varphi(\theta)\,d\nu_y(\theta)\right)
\,d\beta(y).$$
The barycenter of $\beta$ is defined as the barycenter of $\bar\beta$.
\item\label{exlemma3_2} The barycenter is defined in the same way for
  non-negative measures of finite, non-zero mass. 
For any positive constant $c$, we have
  $\bary(c\beta)=\bary(\beta)$.
\item\label{exlemma3_1} 
The barycenter is continuous w.r.t. the weak-$*$ convergence of
  measures, that is, if $\{\beta_i\}$ is a sequence of measures with
  barycenter and converging to a measure $\beta$ with barycenter, then
  $\{\bary(\beta_i)\}\to\bary(\beta)$.
\item\label{exlemma3_3} The barycenter is equivariant by isometries, that is,
  $\bary(\gamma_*\beta)=\gamma(\bary(\beta))$ for any isometry
  $\gamma$ (where $\gamma_*\beta$ denotes the push-forward via $\gamma$ of
  the measure $\beta$).
 \end{enumerate}

\vskip\baselineskip
What we need to complete the proof of Theorem~\ref{mainthm} is the
following fact.  (Compare with Lemma~\ref{mainlemma}.)

\begin{lemma}\label{mainlemma2}
For any $\e>0$ and for any non-elementary representation $\rho$, there exists a map $F^\e:\Hk\to\Hn$ such that
\begin{enumerate}
\item\label{mainlemma2_1}The map $F^\e$ is smooth and $\rho$-equivariant.
\item\label{mainlemma2_2} $|\mathrm{Jac} \,F^\e(x)|\leq 1+\e$, and equality holds if and
  only if $dF^\e_x:T_x\Hk\to T_{F^\e(x)}\Hn$ is a homothety.
\item\label{mainlemma2_3} $\lim_{\e\to0}F^\e=F$, where $F$ is the map of
  Lemma~\ref{mainlemma}.
\item\label{mainlemma2_4} The map $F^\e$ properly ends.
\end{enumerate}
\end{lemma}

Before proving Lemma~\ref{mainlemma2}, we show how it implies
Theorem~\ref{mainthm}.  The inequality directly follows from
points~$(\ref{mainlemma2_1}), (\ref{mainlemma2_2})$ and
$(\ref{mainlemma2_4})$.  If $\vol(\rho)=\vol(M)$, then by point
$(\ref{mainlemma2_3})$ one gets that $\vol(F)=\vol(M)$.  (Note that, a
priori, the map $F$ of Lemma~\ref{mainlemma} does not end
properly).  The proof now follows exactly as in the compact case. 

\begin{proof}[Proof of Lemma~\ref{mainlemma2}.]
The maps $F^\e$ are the so called $\e$-natural maps introduced by
Besson, Courtois, and Gallot.
We begin by recalling their construction.  We omit most details,
referring to~\cite{BCG99,F5,BCG95,BCG96} for a complete discussion on
the construction of natural maps.

For any $\e>0$, we set $$s=(k-1)(1+\e).$$
Let $O$ be a marked point in $\Hk$, and let
$c(s)=\sum_{\gamma\in\Gamma}e^{-sd(O,\gamma O)}$. It turns out that
$c(s)<\infty$, for any $s>k-1$.

For any $x\in\Hk$, we define $\mu_x^\e$ a positive Borel measure on
$\Hk$ by
$$\mu_x^\e=\frac{1}{c(s)}\sum_{\gamma\in\Gamma}e^{-sd(x,\gamma
  O)}\delta_{\gamma O},$$ 
where $\delta_{\gamma O}$ denotes the Dirac measure concentrated on
  the point $\gamma O$.  

Next, we define the measures $\eta_x^\e$ on $\Hn$ and $\lambda_x^\e$  on
$\partial \Hn$, respectively, as the equivariant
push-forward of $\mu_x^\e$ and its convolution with the family
$\{\nu_y\}$ of visual measures. Namely, choose a point $O'\in\Hn$ and
define
$$\eta_x^\e=\frac{1}{c(s)}\sum_{\gamma\in\Gamma}e^{-sd(x,\gamma
  O)}\delta_{\rho(\gamma) O'} \,\,\,\,\,\,\,\,\, \textrm{and}\qquad
\lambda_x^\e=\frac{1}{c(s)}\sum_{\gamma\in\Gamma}e^{-sd(x,\gamma
  O)}\nu_{\rho(\gamma) O'}.
$$ 

The map $F^\e$ is defined by 
$$F^\e(x)=\bary(\eta_x^\e)=\bary(\lambda_x^\e)=
\bary\left(\frac{\lambda_x^\e}{||\lambda_x^\e||}\right).$$ 
Under our present hypotheses we have the following:
\begin{itemize}
\item
(Besson, Courtois, Gallot~\cite[Th\'eor\`eme~1.10]{BCG99}) 
The map $F^\e$ satisfies
  conditions~$(\ref{mainlemma2_1})$ and~$(\ref{mainlemma2_2})$ of
  Lemma~\ref{mainlemma2}. 
\item
(Francaviglia~\cite[Proposition~1.5]{F5}) The maps $F^\e$ satisfy
  condition~$(\ref{mainlemma2_3})$ of Lemma~\ref{mainlemma2}.
\end{itemize}
Therefore, it remains only to prove that for each $\e>0$, the map
$F^\e$ properly ends.
Let $T\subset\overline M$ be a boundary component and let $\pi_1(T)$
be (one of) the corresponding parabolic subgroups of $\pi_1(M)$, and
let $\xi=\mathrm{Fix}(\pi_1(T))$.  

The idea is now the following.
For $x\in\Hk$, we have
$$
\eta_x^\e=\frac{e^{-sd(x,O)}}{c(s)}\sum_{\gamma\in\Gamma}
e^{-s(d(x,\gamma O)-d(x,O))}\delta_{\rho(\gamma)O'},
$$
and by point~$(\ref{exlemma3_2})$ of page~\pageref{exlemma3_2} we have
$$
F^\e(x)=\bary(\eta_x^\e)=\bary\left(\frac{c(s)}{e^{-sd(x,O)}}
\lambda_x^\e\right)=
\bary\left(\sum_{\gamma\in\Gamma}
e^{-s(d(x,\gamma O)-d(x,O))}\nu_{\rho(\gamma)O'}\right).
$$

Now, let $\alpha(t)$ be a geodesic ray ending at $\xi$. As
$t\to\infty$, we have 
$$\sum_{\gamma\in\Gamma}e^{-s(d(\alpha(t),\gamma
  O)-d(\alpha(t),O))}\nu_{\rho(\gamma)O'} 
\stackrel{*}{\rightharpoonup}
\sum_{\gamma\in\Gamma}e^{-sB(\xi,\gamma O)}\nu_{\rho(\gamma)O'},
$$
where $B(\cdot,\cdot)$ denotes the Busemann function normalized at $O$.
Thus, from point~$(\ref{exlemma3_1})$ of page~\pageref{exlemma3_1} we
would get that, as $t\to\infty$
$$F^\e(\alpha(t))\to\bary\left(\sum_{\gamma\in\Gamma}e^{-sB(\xi,\gamma
    O)}\nu_{\rho(\gamma)O'}\right)$$
which one might expect should be fixed by the elements of $\rho(\pi_1(T))$,
because the
    limit measure 
$\sum_{\gamma\in\Gamma}e^{-sB(\xi,\gamma O)}\nu_{\rho(\gamma)O'}$ is
    $\rho(\pi_1(T))$-invariant. 

Unfortunately, the limit measure  
$\sum_{\gamma\in\Gamma}e^{-sB(\xi,\gamma O)}\nu_{\rho(\gamma)O'}$ has
no finite mass, whence its barycenter is not defined.

 In order to overcome this difficulty, some more work is required.
For each $x$ the measure $\lambda_x^\e/||\lambda_x^\e||$ is a
probability measure on $\partial \Hn\simeq \mathbb S^{n-1}$. Since
$\mathbb S^{n-1}$ is compact, the set of probability measures 
on $\partial \Hn$ is weak-$*$ compact.  Therefore, after possibly passing to a subsequence as $x\to\xi$ along the ray $\alpha$, the measures $\lambda_x^\e/||\lambda_x^\e||$ converge to a probability measure $\lambda_\xi$ on $\partial \Hn$.  (The measure $\lambda_\xi$ depends on the chosen subsequence).

We show now that $\lambda_\xi$ is  $\rho(\pi_1(T))$-invariant.
Let $\psi\in\pi_1(T)<\pi_1(M)=\Gamma<\mathrm{Isom}(\Hk)$. Since
$$
\rho(\psi)_*\lambda_x^\e=
\frac{1}{c(s)}\sum_{\gamma\in\Gamma}e^{-sd(x,\gamma
  O)}\nu_{\rho(\psi\gamma) O'}=
\frac{1}{c(s)}\sum_{\gamma\in\Gamma}e^{-sd(x,\psi^{-1}\gamma
  O)}\nu_{\rho(\gamma) O'}
$$
we have
$$
\rho(\psi)_*\lambda_x^\e-\lambda_x^\e=
\frac{1}{c(s)}\sum_{\gamma\in\Gamma}
e^{-sd(x,\gamma O)}
(e^{-s(d(x,\psi^{-1}\gamma O)-d(x,\gamma O))}-1)
\nu_{\rho(\psi\gamma) O'}.
$$
Using the hyperbolic law of sines on the triangles with vertices
$x,\gamma O$ and  $\psi^{-1}\gamma O$, one sees that there exists a function
$E(x)$ such that $E(x)\to 0$ as $x\to \xi$ and
$$|e^{-s(d(x,\psi^{-1}\gamma O)-d(x,\gamma O))}-1|<E(x),$$ 
whence $$||\rho(\psi)_*\lambda_x^\e-\lambda_x^\e||<E(x)||\lambda_x^\e||.$$

Since $||\lambda_x^\e||=||\rho(\psi)_*\lambda_x^\e||$, we have that 
$\lambda_x^\e/||\lambda_x^\e||$ and
$\rho(\psi)_*\lambda_x^\e/||\rho(\psi)_*\lambda_x^\e||$
have the same limit $\lambda_\xi$.
It follows that $\lambda_\xi$ is $\rho(\pi_1(T))$-invariant.

Now we have two cases: either
$\displaystyle{\lambda_\xi=\frac{\delta_{\theta_1}+\delta_{\theta_2}}{2}}$,
or not.
In the latter case, by point~$(\ref{exlemma3_1})$ of
page~\pageref{exlemma3_1},
$$F^\e(x)\to\bary(\lambda_\xi),$$
 which, by 
point~$(\ref{exlemma3_3})$ of page~\pageref{exlemma3_3}, 
is fixed by the elements of $\rho(\pi_1(T))$.  

In the former case, the barycenter of $\lambda_\xi$ is not
defined.  Nevertheless, one can show that the functions
$\B_{\lambda_x^\e}(y)$, defined at page~\pageref{beta}, converge to
$\B_{\lambda_\xi}(y)$.  Since, for each $\e$, $\bary(\lambda_x^\e)$ is the point where $\B_{\lambda_x^\e}$ takes its minimum, they converge to a minimum of $\B_{\lambda_\xi}$ that, by point~$(\ref{newpoint3})$ of page~\pageref{newpoint3}, lies in the geodesic joining $\theta_1$ and $\theta_2$.  Such geodesic is $\rho(\pi_1(T))$-invariant because the invariance of $\lambda_\xi$.  This completes the proof of Lemma
\ref{mainlemma2}, and hence the proof of Theorem \ref{mainthm}.
\end{proof}

\bibliographystyle{plain}

\end{document}